\documentclass[letterpaper, 10 pt, conference]{ieeeconf}  

\IEEEoverridecommandlockouts                              

\usepackage{graphics} 
\usepackage{epsfig} 
\usepackage{mathptmx} 
\usepackage{times} 
\usepackage{amsmath} 
\usepackage{amssymb}  
\usepackage{svg}
\usepackage{subfigure}
\usepackage{cite}
\usepackage{mathrsfs}
\newcommand{\myproof}[1][Proof]{%
    \par\noindent\hspace{2em}\textit{\textbf{#1: }}\ignorespaces}
\newcommand{\myqed}{\hfill $\square$}

\title{\LARGE \bf
Differential Game with Geometric Formulation for Differential Systems Under Bilateral Perception Constraints\thanks{This work has been submitted to the IEEE for possible publication. Copyright may be transferred without notice, after which this version may no longer be accessible.}}

\author{Xinyi Zhu, Jiali Wang, Yang Tang, \IEEEmembership{Fellow, IEEE}, Fangfei Li, \IEEEmembership{Member, IEEE}, and Yan Zhu
\thanks{This work is supported by the National Natural Science Foundation of China (No. 62233005, 62293502, 62173142, U2441245), the Programme of Introducing Talents of Discipline to Universities (the 111 Project) (No.B17017), and the National Key Laboratory of Space Intelligent Control (No. HTKJ2024KL502004).}
\thanks{X. Zhu and Y. Zhu are with the School of Mathematics, East China University of Science and Technology, 130 Meilong Road, Shanghai 200237, China (e-mail: zxyihemmm@163.com, zhuygraph@ecust.edu.cn) }
\thanks{J. Wang and Y. Tang are with Key Laboratory of Smart Manufacturing in Energy Chemical Process, Ministry of Education, East China University of Science and Technology, Shanghai 200237, China (e-mail: jialiwang@mail.ecust.edu.cn, yangtang@ecust.edu.cn).}
\thanks{F. Li is with the School of Mathematics and the Key Laboratory of Smart Manufacturing in Energy Chemical Process, Ministry of Education, East China University of Science and Technology, Shanghai 200237, China (email: li\_fangfei@163.com)}}

\begin{document}

\maketitle
\thispagestyle{empty}
\pagestyle{empty}

\begin{abstract}

This letter employs differential game theory to address the defense problem of a circular target area with perception constraints, involving a single defender and a single attacker. The defender is restricted to moving along the perimeter, while the mobile attacker aims to make itself to the edge of the circular target to win. We examine a scenario where both the attacker and defender face perception constraints, dividing the interaction into four distinct stages based on detection capabilities and deriving the corresponding optimal control strategies. Simulations are conducted to validate the proposed strategies.

\end{abstract}


\section{INTRODUCTION}

Evasive targets, which require detection, interception, or containment by others, are modeled as pursuit-evasion games \cite{r5}. Since Isaacs' pioneering work on differential games \cite{r6}, these games have been widely studied \cite{r8,r10,r11,r12}, with one such application being the target defense problem. In the target defense problem, a classic in differential game theory, a defender and an attacker contest a target. The defender's goal is to block the attacker, who seeks to reach the target using optimal strategies. Differential games provide a dynamic framework for modeling multi-agent systems \cite{r27}, with strategies represented by state-space trajectories. Target defense problems have applications in military defense \cite{r1}, aerospace \cite{r2}, and cybersecurity \cite{r4}.

Various target defense problems have been studied using different methods \cite{r13,r14,r15,r16,r17}. For example, \cite{r19} applied the Hamilton-Jacobi-Isaacs (HJI) solution and path defense to one-on-one Target-Attacker-Defender (TAD) problems, while \cite{r20} used deep reinforcement learning. Circular target defense problems, in which both the attacker and defender can detect each other, were explored in \cite{r21}.

Although the study in \cite{r21} explored an interesting issue, we note that the literature in \cite{r21} considers a full-information game and does not account for perception constraints. In fact, perception constraints are common in real life, and incomplete information complicates game analysis \cite{r25}. As noted in \cite{r24}, when an aircraft enters an unknown or unmonitored space, the collision risk increases significantly. Thus, incorporating perception constraints is crucial for realistic modeling in TAD games.

Incorporating perception constraints into TAD games adds complexity, as nonlinear differential equations become harder to solve, and numerical methods may struggle with stability. Different constraints affect game dynamics in varying ways, and decision-making under incomplete information leads to large-scale dynamic programming challenges. Multi-stage games with perception constraints also require careful analysis to ensure system continuity and strategy coherence.

Building on \cite{r21}, we introduce perception constraints for both factions and analyze scenarios with limited perception. The process is divided into three stages: the Pre-game Stage, where the attacker approaches the Target Sensing Region (TSR) \cite{r5} without detection; the Partial-information Stage, where the attacker enters the TSR but cannot perceive the defender; and the Full-information Stage, where both sides engage in a Full-information Game with derived equilibrium state feedback strategies. If the defender wins, an Escape Stage is added to evade the attacker.

The contribution of this letter is threefold. First, we incorporate perception constraints for both factions, refining optimal strategies and reflecting real-world scenarios, unlike the complete information game in \cite{r21}. Secondly, compared to \cite{r13}, we introduce a geometric perspective, which offers a clearer geometric interpretation and a global view for problem analysis. This geometric approach effectively captures the multi-stage dynamics induced by perception constraints, offering a clearer representation of critical system features such as the shape of optimal strategy trajectories and the geometric structure of winning regions for both factions, which are less explicitly addressed in the method of \cite{r13}. Third, compared with \cite{r21}, we conduct numerical simulations of the game with perception constraints, and provide a comparative analysis of how these constraints impact full-information game outcomes.

The remainder of this letter is organized as follows: Section II describes the problem and outlines the assumptions. Section III discusses the scenario where both the defender and attacker have limited perception. Section IV presents the simulation and emulation of the strategic derivations. Section V summarizes and concludes the letter.

\textbf{Notations:} Variables are nondimensionalized for clarity. The system state is denoted by $\overline{\mathbf{s}}$, with its nondimensionalized form as $\mathbf{s}$. The coordinates of the attacker and defender are $\mathbf{x}_A$ and $\mathbf{x}_D$, respectively. Time is represented by $\overline{t}$ and is nondimensionalized as $t$. The distance from the attacker to the target center is $\overline{R}$, with its nondimensionalized form as $R$. The angle between the attacker and defender is $\overline{\theta}$, nondimensionalized as $\theta$. The defender’s rotation around the circle’s center is $\overline{\beta}$, with its nondimensionalized counterpart denoted by $\beta$. The attacker’s control is $\alpha$, and the defender’s control, after nondimensionalization, is $u$. At termination, the nondimensionalized distance from the attacker to the target is $R_f$, and the system state is $\mathbf{s}_f$. The termination time is $t_f$, and the terminal separation angle is $\theta_f$. The outer radius of TSR is $\overline{R}_0$, with its nondimensionalized version as $R_0$. The initial nondimensionalized distance from the attacker to the target center is $L_0$. Let $\theta_p$ denote the angle of the attacker's position to the defender when the defender reaches the edge of the ASR, with the variables nondimensionalized.
\section{Problem Description}
This letter addresses the defense of a circular target under bilateral perception constraints. The target is a circle with a radius $L$. The defender, $\mathcal{D}$, patrols the perimeter to intercept the attacker, $\mathcal{A}$, who attempts to breach the target and reach its circumference. Perception actions enable the detection of state information, including positions and velocities, within sensing areas. The defender's sensing range is a ring around the target, with an outer radius $\overline{R}_0$ and an inner radius $L$. The attacker's detection range (ASR) is a circle centered on the attacker, with radius $\overline{r}_A$. After nondimensionalization, $r_A$ is defined as $r_A \equiv \frac{\overline{r}_A}{L}.$

$\overline{v}_1$ represents the attacker's velocity and $\overline{v}_2$ is the defender's velocity. If $\overline{v}_1 > \overline{v}_2$, the attacker only needs to get within a certain distance where they can control the angle, forcing $\theta \rightarrow \pi$ \cite{r21}, meaning the defender cannot close the angle gap. As a result, interception at the target's edge becomes impossible. Therefore, we make the following assumptions:

Assumption 1 \cite{r21}: The velocities satisfy $0\leq\overline{v}_1\leq \overline{v}_2$. 

Assumption 2 \cite{r21}: The initial separation angle satisfies $\overline \theta(t_{0})=\overline \theta_{0}\in [0,\pi)$.

Assumption 3: The initial distance between the attacker and the target center is $\overline{R}(t_0) \geq L$, meaning $\mathcal{A}$ starts outside the circular target.

To simplify the analysis and standardize the scale, the relative variables are nondimensionalized using the following definitions:
\begin{equation}
    R \equiv \frac{\overline{R}}{L},t \equiv \frac{\overline{v}_{2max}}{L}\overline{t},u \equiv \frac{\overline{v}_{2}}{\overline{v}_{2max}},\nu \equiv \frac{\overline{v}_{1}}{\overline{v}_{2max}},
    \theta \equiv \overline{\theta},
    \beta \equiv \overline{\beta},
    \nonumber
\end{equation}
where $\overline{v}_{2max}$ represents the defender's maximum possible speed, while the speed ratio $\nu$ satisfies $\nu \in (0,1]$. The defender's control is $u \in [-1, 1]$, and the attacker's control is $\alpha \in [-\pi, \pi]$. Both controls are independent of time and constrained to their respective ranges. Since the attacker aims to approach the target, they will inevitably enter the TSR. Thus, we can further restrict $\alpha$ to the range $(-\frac{\pi}{2}, \frac{\pi}{2})$.

Assumption 4: In real-world applications, the target typically covers a large area, while the attacker is a compact, mobile entity. Thus, the attacker's perception radius is smaller than the distance between the D-T camp perimeter and the target's perimeter, i.e., $0<r_A\leq\min{\left\{R_0-1,1\right\}}$.

\begin{figure}[h]
\centerline{\includegraphics[width=\columnwidth]{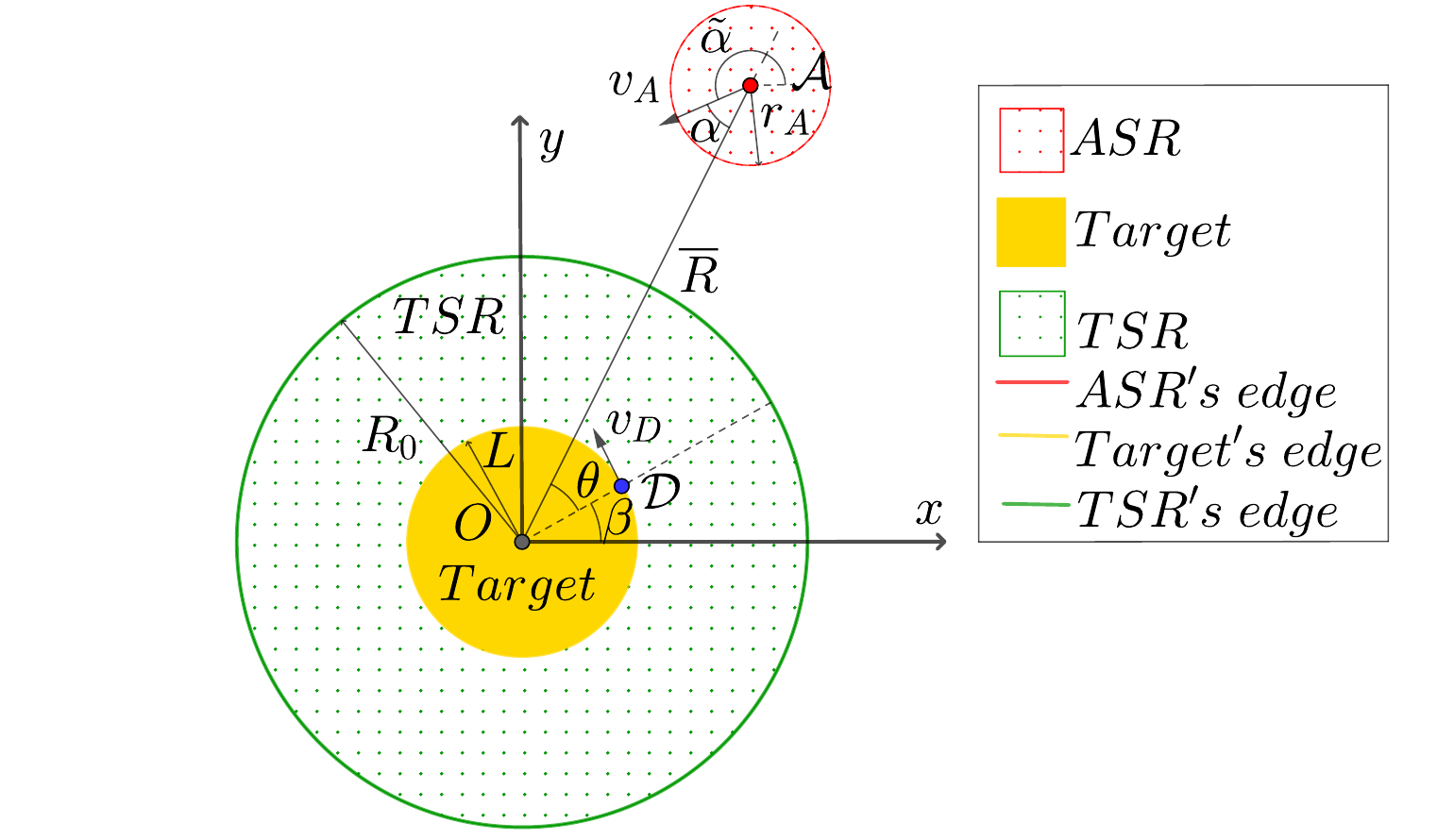}}
\caption{\centering{Circular perimeter patrol when both attacker and defender have perception constraints.}}
\label{fig1}
\end{figure}
\begin{figure}[h]
\centerline{\includegraphics[width=\columnwidth]{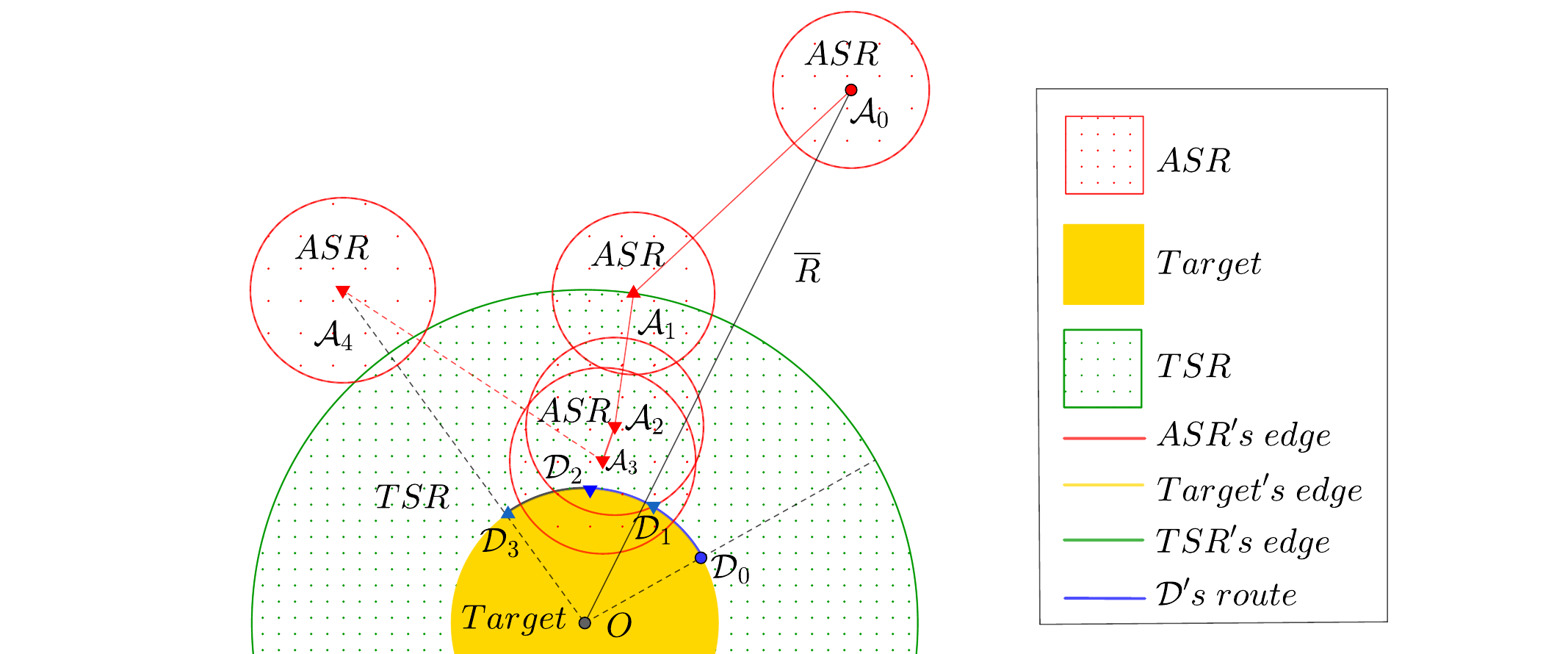}}
\caption{\centering{The entire process of TAD games with limited players' vision.}}
\label{fig1}
\end{figure}

\section{Perception Constraints on Both the D-T Factions and the Attacker}
As shown in Fig. 1, the defender's sensory capabilities are limited to the TSR, meaning that the attacker can only be detected by the defender upon entering the TSR. Meanwhile, their perception is confined to the attacker's sensing region (ASR), which implies that the defender is only perceived by the attacker when within the ASR. Under these perception constraints, the game progresses through the stages outlined in Fig. 2.

\subsection{The Pre-game Stage}
In this stage, the attacker has not entered the TSR, so the defender cannot perceive them, and remains still. During this phase, the attacker's starting position and path to the TSR are entirely unpredictable. The details are shown in Fig. 2, where the attacker moves from $\mathcal{A}_0$ to $\mathcal{A}_1$ and the defender remains $\mathcal{D}_0$. The defender stays stationary while the attacker advances freely toward the TSR. The initial distance from the attacker to the target center is $R(t_0) = R_0$, which remains unchanged at the end of the Pre-game Stage and the start of the Full-information Game. 

\subsection{The Partial-information Stage}

The partial-information stage begins when the attacker approaches the TSR edge and ends when the defender reaches the ASR edge. During this stage, the defender can perceive the attacker, while the attacker remains unaware of the defender. As shown in Fig. 2, during the Partial-information Stage, the attacker moves from $\mathcal{A}_1$ to $\mathcal{A}_2$, and the defender from $\mathcal{D}_0$ to $\mathcal{D}_1$. The attacker, unaware of the defender's exact location, proceeds without adjusting for the defender's actions. 
The dynamics are as follows:
\begin{equation}
f=\left(\mathbf{s},u,t\right)=\dot{\mathbf{s}}=\begin{bmatrix}\dot{R}(t)\\\dot{\theta}(t)\end{bmatrix}=\begin{bmatrix}-\nu\cos{\alpha}\\\frac{\nu}{R(t)}\sin{\alpha}\end{bmatrix}.
\label{eq} 
\end{equation}
To make $R \rightarrow 1$ as quickly as possible, the cost function of the attacker should be
\begin{equation}
    J_a=\Phi_a(\mathbf{s}_f,t_f)=-t_f.
    \label{eq}
\end{equation}
This is called the \textit{Game of Time}, symbolized by the subscript `$a$'. The termination criterion for this game is given by 
\begin{equation}
    \phi_a(\mathbf{s}_f,t_f)=R_f-1=0.
    \label{eq}
\end{equation}
The condition for the Partial-information Stage is:
\begin{equation}
    \theta>\theta_p,
    \label{eq}
\end{equation}
where $\theta_p>0$. 

\textit{Theorem 1}: In the differential game with the cost defined by (2), the equilibrium control of the defender is
\begin{equation}
    u^{*}=\operatorname{sign}(\theta),
    \label{eq}
\end{equation}
while the attacker's optimal control strategy is move towards the target along the line connecting their current position to the center of the target.

\myproof
    The defender's main goal is to intercept the attacker quickly, driving $\theta \rightarrow 0$. The most efficient way is to advance directly toward the attacker at maximum speed, as shown in equation (5). As for the attacker, the Hamiltonian is expressed as
\begin{equation}
\mathscr{H}_{a}=-\sigma_{R}\nu \cos{\alpha}+\sigma_{\theta}\frac{\nu}{R} \sin{\alpha},
\label{eq} 
\end{equation}
where $\sigma \equiv \begin{bmatrix} \sigma_R &\sigma_\theta\end{bmatrix}^\top$. Using the method in \cite{r21}, one has:
\begin{equation}
    \alpha=\arccos({-\operatorname{sign}(\sigma_R)})=0\;or\; \pi.
    \label{eq}
\end{equation}  
During the Partial-information Stage, the attacker should move towards the target along the line connecting their position to the target's center.
\myqed

\subsection{The Full-information Stage}

In the Full-information Stage, both the attacker in the TSR and the defender in the ASR can perceive each other, making it a Full-information Game (as shown in Fig. 2, where the attacker moves from $\mathcal{A}_2$ to $\mathcal{A}_3$, and the defender from $\mathcal{D}_1$ to $\mathcal{D}_2$). The dynamics are \cite{r21}:
\begin{equation}
f=\left(\mathbf{s},u,t\right)=\dot{\mathbf{s}}=\begin{bmatrix}\dot{R}(t)\\\dot{\theta}(t)\\\dot{\beta}(t)\end{bmatrix}=\begin{bmatrix}-\nu\cos{\alpha}\\\frac{\nu}{R(t)}\sin{\alpha}-u\\u\end{bmatrix}.
\label{eq} 
\end{equation}

However, this stage may not always occur. If the defender starts far from the attacker, the attacker might not detect the defender even upon reaching the target's perimeter. Given an initial distance of $R(0) = L_0$ for the attacker, and assuming the Full-information Stage is necessary, the defender should have entered the ASR by the time the attacker reaches the target's edge, following the Partial-information Stage strategy. Let the initial angular difference be $\theta_0$. During the Partial-information Stage, the attacker moves along the line toward the target's center, implying $\alpha = 0$ and $\dot{\theta} = u = -1$.

\textit{Theorem 2}: To ensure the presence of the Full-information Stage, $\theta_0$ must satisfy the following condition:
$
\theta_0 \in \Big[-\frac{{L_0 - 1}}{\nu}-\arccos \left( 1 - \frac{{r_A^2}}{2} \right),\frac{{L_0 - 1}}{\nu}+\arccos \left( 1 - \frac{{r_A^2}}{2} \right)\Big].
$

\myproof
In the critical scenario where the attacker reaches the target's edge, and the defender just enters the attacker's perceptual range (ASR), the time taken by the attacker to reach the terminal point is $ t = \frac{{L_0 - 1}}{\nu} $. Thus, the defender's operating angle is $\theta = \frac{{v_2 t}}{L} = u t = t $. Upon reaching this critical state, the heading difference angle is given by the cosine rule: $\theta_1 = \arccos \left( 1 - \frac{{r_A^2}}{2} \right)$.
Therefore, the condition for the existence of the Full-information Stage is $\theta_0 - t = \theta_1$, which means
\begin{equation}
    \theta_0=\frac{{L_0 - 1}}{\nu}+\arccos \left( 1 - \frac{{r_A^2}}{2} \right).
    \label{eq}
\end{equation}
Due to symmetry, an additional critical scenario aries as
\begin{equation}
    \theta_0=-\frac{{L_0 - 1}}{\nu}-\arccos \left( 1 - \frac{{r_A^2}}{2} \right).
    \label{eq}
\end{equation}
The condition for the existence of the Full-information Stage is that the initial angular difference, $\theta_0$, satisfies
\begin{small}
\begin{equation}
    \theta_0 \in \big[-\frac{{L_0 - 1}}{\nu}-\arccos \left( 1 - \frac{{r_A^2}}{2} \right),\frac{{L_0 - 1}}{\nu}+\arccos \left( 1 - \frac{{r_A^2}}{2} \right)\big].
    \label{eq}
\end{equation}
\end{small}
\myqed

The distance between the attacker and defender, denoted as $\overline{p}$, can be derived using the cosine theorem, resulting in 
\begin{equation}
    \overline{p}=\sqrt{\overline{R}^2+L^2-2\overline{R}L\cos{\overline{\theta}}}.
    \label{eq}
\end{equation}
Upon nondimensionalizing $p$ such that $p\equiv\frac{\overline{p}}{L}$, we obtain  
\begin{equation}
    p=\sqrt{R^2+1-2R\cos{\theta}}.
    \label{eq}
\end{equation}
To find the value of $\theta$ when the defender enters the attacker's perception range, let $p = r_A$. This angle is denoted as $\theta_p$, leading to the result:
\begin{equation}
    \theta_p=\arccos{\frac{R^2+1-r^2_A}{2R}},R\in[1,1+r_A].
    \label{eq}
\end{equation}

\textit{Proposition 1}: The range of $\theta_p$ is $
\theta_p\in (0,\arccos{(\frac{2-r^2_A}{2})}]$.
\myproof
Given that
$$\frac{R^2+1-r^2_A}{2R}\in[\frac{2-r^2_A}{2},1]\subset (\frac{1}{2},1],$$ it follows that,$$\theta_p\in [0,\arccos{(\frac{2-r^2_A}{2})}].$$
Considering the corner case $\theta=0$, we have,
\begin{equation}
    \frac{R^2+1-r^2_A}{2R}=1.
    \label{eq}
\end{equation}
Solving (15), yields
\begin{equation}
    R=1 \pm r_A.
    \label{eq}
\end{equation}
Since $R\geq 1$, we conclude that
\begin{equation}
    R=1+r_{A}.
    \label{eq}
\end{equation}
At the end of the Partial-information Stage, the dynamics of $R$ and $\theta$ are governed by
\begin{equation}
    \dot{R}=-\nu,\dot{\theta}=-1.
    \label{eq}
\end{equation}
This implies that $\frac{d\theta}{dR} = \frac{1}{\nu}$, meaning the trajectories in the $(R, \theta)$ plane are straight lines. Additionally, the unconstrained equilibrium flow field is described by
\begin{equation}
    \theta(R)=\frac{1}{\nu}(R-R_0)+\theta_0.
    \label{eq}
\end{equation}
The critical condition occurs when the attacker reaches the target ($\theta_f = 0$), activating the constraint. The time taken from their initial locations must be the same. By deducing backwards from the end of the Partial-information Stage, we obtain 
\begin{equation}
    \theta_c(R)=\frac{R-1}{\nu}.
    \label{eq}
\end{equation}
Substituting $R=1+r_A$ into (20), we get
\begin{equation}
    \theta_c(1+r_A)=\frac{r_A}{\nu}>0.
    \label{eq}
\end{equation}
Therefore, let $R = r_A + 1$, and $\theta > 0$, indicating that $\theta_p$ cannot reach 0. Hence,
\begin{equation}
    \theta_p\in (0,\arccos{(\frac{2-r^2_A}{2})}].
\end{equation}
\myqed

Considering (8), where $\dot{\theta} < 0$, we see that $\theta$ decreases monotonically. Once $\theta$ reaches $\theta_p$, the defender becomes continuously detectable by the attacker. Thus, when $\theta > \theta_p$, the attacker cannot detect the defender; when $\theta \leq \theta_p$, the attacker can sense the defender. After the Partial-information Stage transitions into the Full-information Stage, the latter will continue indefinitely.

\textit{Theorem 3}: If the Full-information Stage exists, the strategy of the attacker is:
\begin{equation}
\begin{cases}
\alpha^{*}=\operatorname{sign}(\theta)\arcsin{(\frac{\nu}{R})}, & \text{if } \theta_G\leq \theta\leq\theta_p \\
\text{stop attacking}, & \text{if } \theta<\theta_p
\end{cases}
    \label{eq}
\end{equation}
while the defender’s optimal strategy is the same as (5).

\myproof
    When $\theta < \theta_p$, the Full-information Stage begins. As outlined in \cite{r21}, the equilibrium state feedback control strategies, when both the defender and attacker can perceive each other(i.e., $\theta \leq \theta_p$), are
\begin{equation}
    \alpha^{*}=\operatorname{sign}(\theta)\arcsin{(\frac{\nu}{R})}\; \textup{and} \;u^{*}=\operatorname{sign}(\theta),
    \label{eq}
\end{equation}
and the attacker follows a straight-line path in the inertial $(x, y)$-plane. The surface  
\begin{equation}
    \theta_{G}(R)=g(R)-g(1),\forall R\in (1,\frac{R_0}{L}],
    \label{eq}
\end{equation}
where $g(R)=\sqrt{\frac{R^2}{\nu^2}-1}+\arcsin{(\frac{\nu}{R})}$,
divides the state space into regions of win for the defender and attacker:
\begin{equation}
    \mathcal{R}_{D}=\left\{\mathbf{s}||\theta|\leq\theta_{G}(R)\right\},
\end{equation}
\begin{equation}
    \mathcal{R}_{A}=\left\{\mathbf{s}||\theta|>\theta_{G}(R)\right\}.
\end{equation}
Combining with $\theta_G$ and \textit{Theorem 1}, if the system state falls within $\mathcal{R}_{D}$, the defender is certain to win, and to minimize losses, the attacker should cease the attack to preserve itself. On the other hand, if the system state enters $\mathcal{R}_{A}$, the attacker is assured victory and should therefore continue the attack.
\myqed

\subsection{The Escape Stage}
If the Full-information Stage occurs and the attacker fails to penetrate, they will withdraw beyond the TSR. The defender will maintain $\theta = 0$ to prevent a reversal (as shown in Fig. 2, with the attacker moving from $\mathcal{A}_2$ to $\mathcal{A}_3$ and the defender from $\mathcal{D}_1$ to $\mathcal{D}_2$). The dynamics resemble those in equation (8), where $\alpha \in [\frac{\pi}{2}, \frac{3\pi}{2}]$. Since the attacker's objective is to reach the outer boundary of the TSR promptly, the cost functional is
\begin{equation}
    J_E=\Phi_{E}(x_{f},t_{f})=-t_f,
    \label{eq}
\end{equation}
where the subscript `E' represents the escape stage. The attacker seeks to maximize (28), while the defender aims to minimize it, following the previously established conventions. Since the defender aims to maintain $\theta = 0$, it follows that $\dot{\theta}(t) = 0$. Therefore, the defender must ensure that $u$ satisfies the following condition:
\begin{equation}
    u=\frac{\nu}{R}\sin{\alpha}.
    \label{eq}
\end{equation}
The terminal constraint is
\begin{equation}
\phi_{E}(x_{f},t_{f})=R - \frac{R_0}{L}=0.
\label{eq} 
\end{equation}
The earliest moment when $R(t) = \frac{R_0}{L}$ is referred to as the final time $t_f$. Thus, the terminal surface corresponds to the set of states that satisfy equation (31):
\begin{equation}
\mathcal{J}_{E}=\{ x \mid R = \frac{R_0}{L} \text{ and } \theta = 0 \}.
\label{eq} 
\end{equation}

\begin{figure*}
    \begin{minipage}[t]{0.31\linewidth}
        \centering
        \subfigure[]{
            \includegraphics[width=1\linewidth]{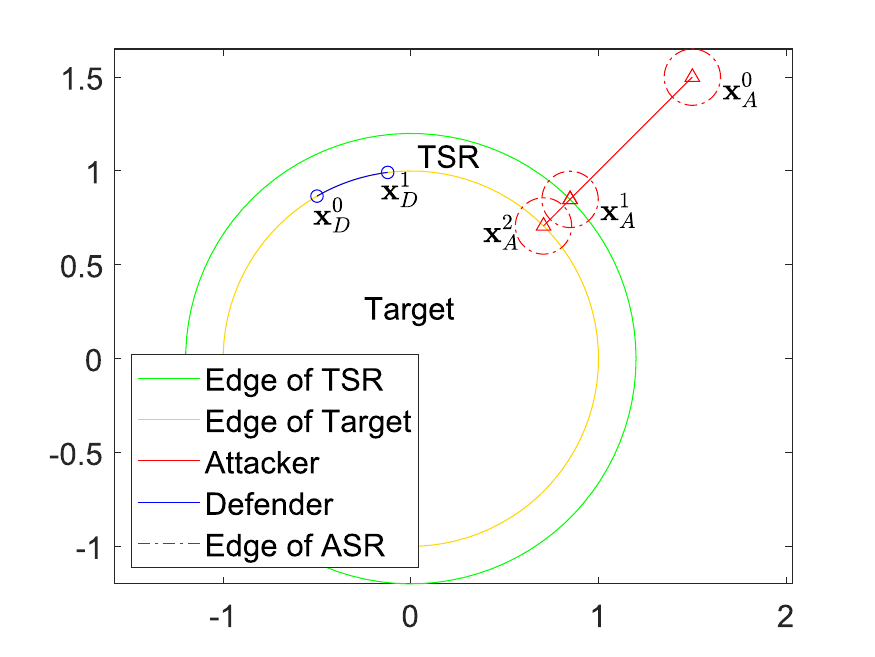}
        }

        \subfigure[]{
            \includegraphics[width=1\linewidth]{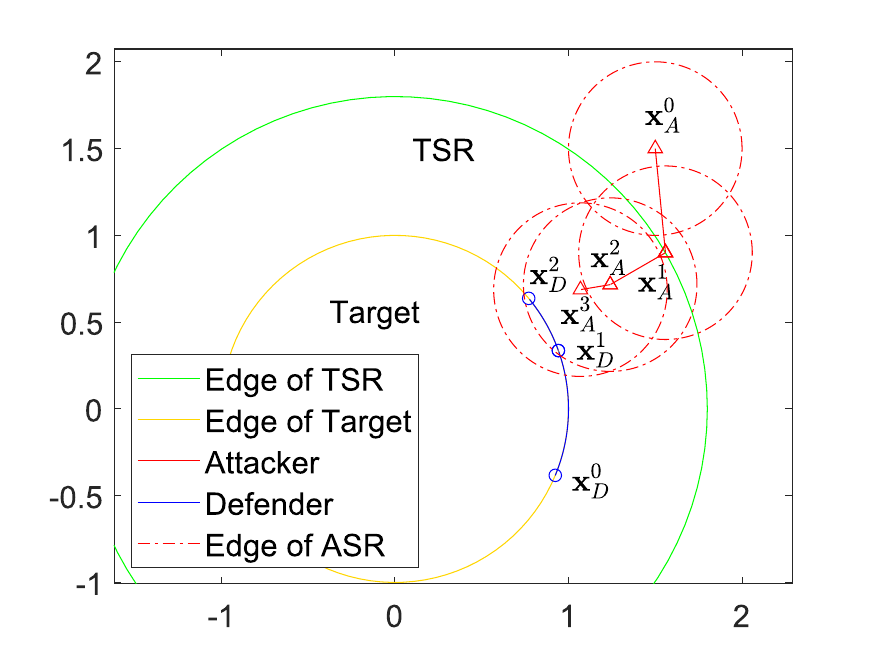}
        }
        \caption{Both the attacker and defender face perception constraints.}
    \end{minipage}
    \hfill
    \begin{minipage}[t]{0.31\linewidth}
        \centering
        \subfigure[]{
            \includegraphics[width=1\linewidth]{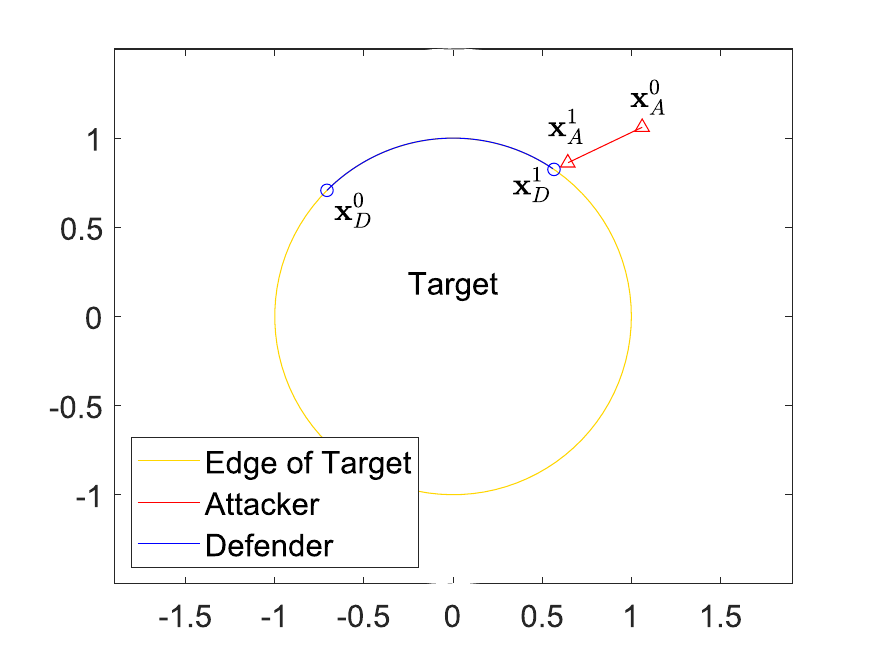}
        }
        
        \subfigure[]{
            \includegraphics[width=1\linewidth]{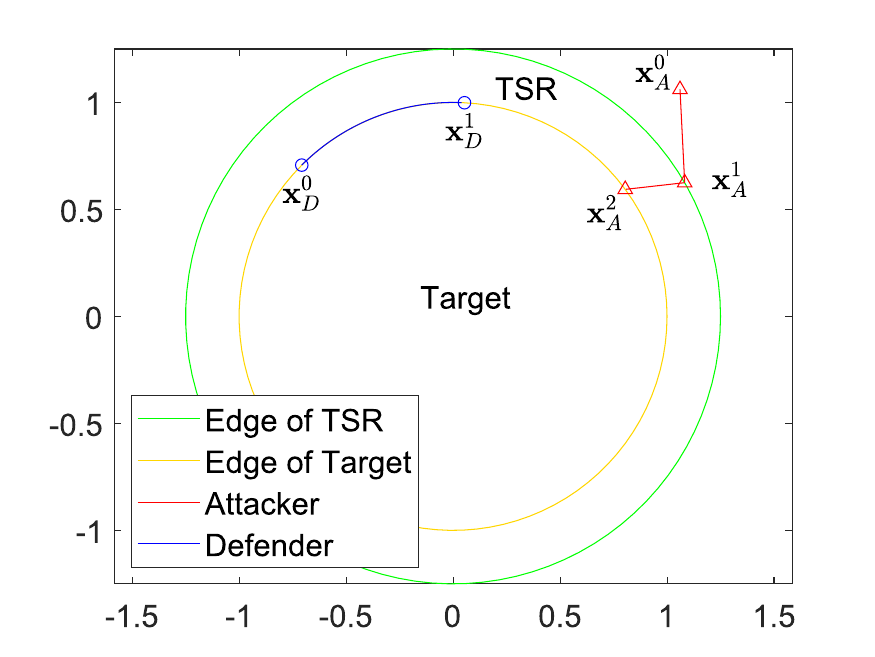}
        }
        \caption{The impact of adding TSR.}
    \end{minipage}
    \hfill
    \begin{minipage}[t]{0.31\linewidth}
        \centering
        \subfigure[]{
            \includegraphics[width=1\linewidth]{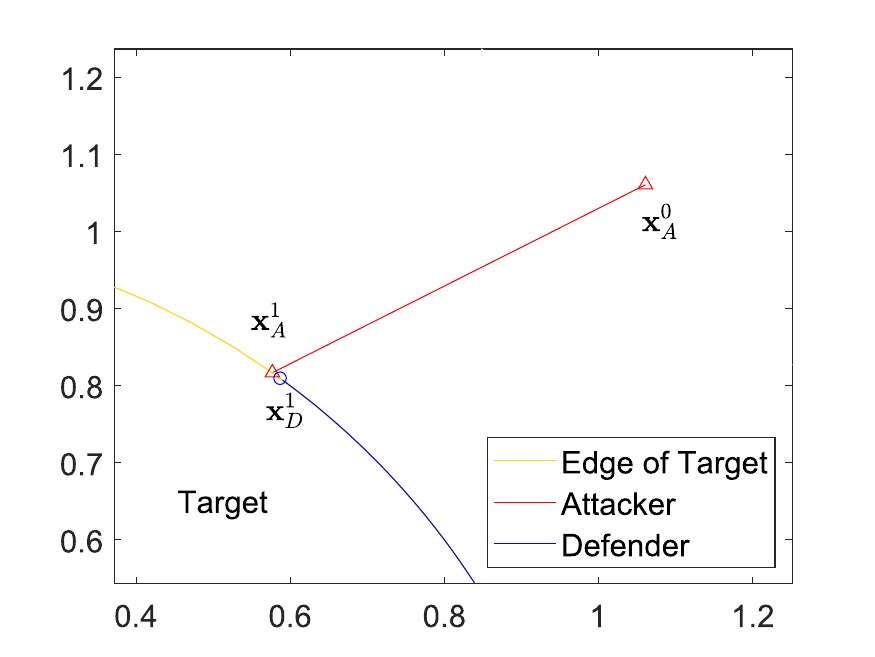}
        }
        
        \subfigure[]{
            \includegraphics[width=1\linewidth]{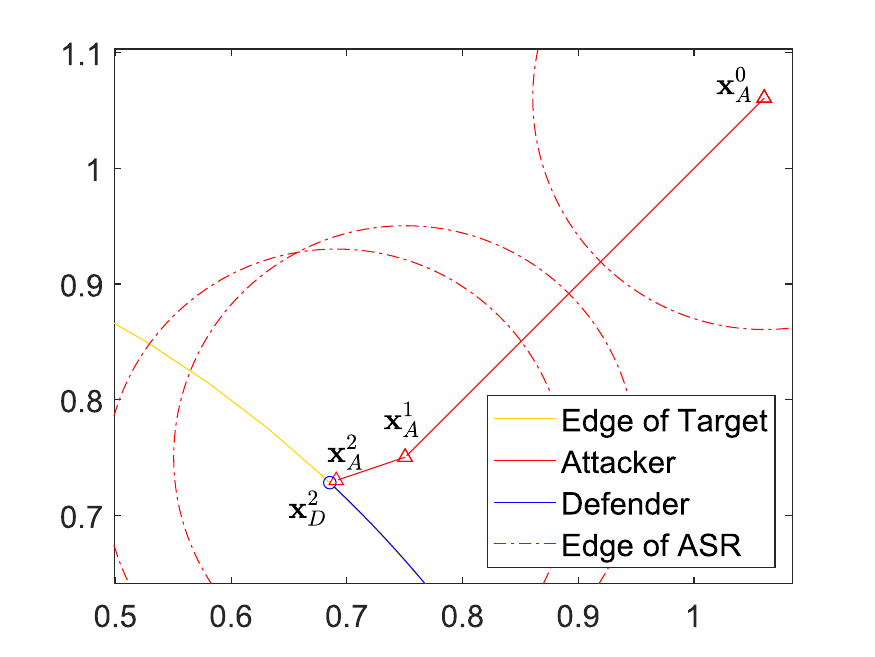}
        }
        \caption{The impact of adding ASR.}
    \end{minipage}
\end{figure*}

\textit{Theorem 4}: The optimal escape strategy for the attacker is to move outward along the line connecting its position to the center of the target, which is represented by $\alpha^* = \pi$. Meanwhile, the defender should remain stationary, which is indicated by $u = 0$.

\myproof
The dynamics are given by equation (8), and the Hamiltonian is:
\begin{equation}
\mathscr{H}_{E}=-\mu_{R}\nu \cos{\alpha}+\mu_{\theta}(\frac{\nu}{R} \sin{\alpha}- u)+\mu_{\beta} u.
\label{eq} 
\end{equation}
Since (29), we can obtain
\begin{equation}
    \mathscr{H}_{E}=-\mu_{R}\nu \cos{\alpha}+\mu_{\beta} \frac{\nu}{R}\sin{\alpha},
    \label{eq}
\end{equation}
where $\mu \equiv \begin{bmatrix} \mu_R & \mu_\beta \end{bmatrix}^\top$ is the adjoint vector. And the corresponding equilibrium adjoint dynamics are given by
\begin{equation}
\dot{\mu}_{R} = - \frac{\partial \mathscr{H}_{E}}{\partial R} = \frac{\nu\mu_\theta}{R^{2}} \sin{\alpha},
\label{eq} 
\end{equation}
\begin{equation}
\dot{\mu}_{\beta} = - \frac{\partial \mathscr{H}_{E}}{\partial \beta} = 0.
\label{eq} 
\end{equation}
The terminal adjoint values are determined using the transversality condition \cite{r23}, which is given as follows
\begin{equation}
\mu^\top(t_{f})=-\frac{\partial \Phi_{E}}{\partial\mathbf{s}_f} +\sigma \frac{\partial \phi_{E}}{\partial\mathbf{s}_f}\\
=\begin{bmatrix}\sigma &  0 \\\end{bmatrix}.
\end{equation}
Based on (36), we have
\begin{equation}
\mu_{R_f} =\sigma, 
\mu_{\beta_f} = 0,
\label{eq} 
\end{equation}
where $\sigma$ is an additional adjoint variable. Therefore, from (35) and (37), the following holds
\begin{equation}
\mu_{\beta} = 0,\forall t \in [t_0, t_f].
\label{eq} 
\end{equation}
According to (38), the state component $\beta$ does not affect the equilibrium trajectory or control strategies. We can rewrite (39) as:
\begin{equation}
\mathscr{H}_{E}=-\mu_{R}\nu \cos{\alpha}.
\label{eq} 
\end{equation}
The attacker aims to maximize (39), therefore
\begin{equation}
\alpha^*=\pi,u^*=\frac{\nu}{R}\sin{\alpha}=0,
\label{eq} 
\end{equation}
which implies that the defender remains stationary. Furthermore, $\theta = 0$ holds true throughout the escape process.
\myqed

\section{ILLUSTRATIVE EXAMPLES}
Next, we conduct simulations and modeling on a two-dimensional plane, where both the $x$-axis and $y$-axis represent dimensionless lengths. Actual distances are normalized by dividing by $L$, the target area's radius, simplifying the model by converting lengths into unitless quantities.

\subsection{The Feasibility of Strategies}
We simulate the outcomes as follows:

1) The attacker wins when both the attacker and defender have limited perception, and no third stage occurs. In Fig. 3 (a), the radius of the TSR is $1.2$, the radius of the target is 1, and the radius of the ASR is 0.15. The speed ratio is $\nu = 0.5$. The initial positions are: $\mathcal{D}$ starts at $\mathbf{x}_D^0 \left( -\frac{1}{2}, \frac{\sqrt{3}}{2} \right)$, and $\mathcal{A}$ starts at $\mathbf{x}_A^0 \left( 1.5, 1.5 \right)$. $\mathcal{A}$ reaches the edge of TSR at $\mathbf{x}_A^1 \left( \frac{3\sqrt{2}}{5}, \frac{3\sqrt{2}}{5} \right)$. The position where $\mathcal{A}$ reaches the target's edge is denoted by $\mathbf{x}_A^2$. At this moment, $\mathcal{D}$ is at $\mathbf{x}_D^1$ and fails to intercept $\mathcal{A}$, so $\mathcal{A}$ wins.

2) The scenario includes an escape stage when both the attacker and defender have limited perception and a third stage is present. In Fig. 3 (b), the TSR radius is 1.8, the target radius is 1, and the ASR radius is 0.5. The speed ratio is $\nu = 0.5$. The initial position of $\mathcal{D}$ is $\mathbf{x}_D^0 \left( \cos{(-\frac{\pi}{8})}, \sin{(-\frac{\pi}{8})} \right)$, and $\mathcal{A}$ starts at $\mathbf{x}_A^0 \left( 1.5, 1.5 \right)$. $\mathcal{A}$ reaches the edge of the TSR at $\mathbf{x}_A^1 \left( \frac{9\sqrt{3}}{10}, \frac{9}{10} \right)$.
Special points: (1) When $\mathcal{D}$ enters $\mathcal{A}$'s perception at $\mathbf{x}_D^1$ and $\mathcal{A}$ is at $\mathbf{x}_A^2$, $\mathcal{A}$ switches to a full-information strategy. $\mathcal{D}$ continues along the edge of T toward $\mathcal{A}$. (2) At $\mathbf{x}_A^3$, $\mathcal{A}$ realizes that $\mathcal{D}$ has reached the endpoint of its path at $\mathbf{x}_D^2 (\theta = 0)$. Consequently, $\mathcal{A}$ adopts an escape strategy and moves to $\mathbf{x}_A^4$ outside the TSR. Throughout, $\mathcal{D}$ maintains a constant heading difference of $\theta = 0$ relative to $\mathcal{A}$, remaining stationary.

\subsection{The Impact of Adding Perception Constraints}
The inclusion of perception constraints significantly affects the final game outcome. In this subsection, we select several scenarios from reference \cite{r21} and compare them with the results obtained after incorporating perception constraints, including TSR and ASR.

1) The addition of TSR: The initial states are as follows: $\mathcal{D}$ is at $\mathbf{x}_D^0 \left( -\frac{1}{\sqrt{2}}, \frac{1}{\sqrt{2}} \right)$ and $\mathcal{A}$ is at $\mathbf{x}_A^0 \left( \frac{3}{2\sqrt{2}}, \frac{3}{2\sqrt{2}} \right)$. The speed ratio is $\nu = 0.5$, and the target radius is 1. In Fig. 4 (a), both sides use the strategy from \cite{r21}, and the defender successfully intercepts the attacker. However, when a defender perception radius of 1.25 (TSR) is introduced, as shown in Fig. 4 (b), the interception fails. This indicates that the successful interception strategy in \cite{r21} becomes ineffective with the addition of TSR. The strategy that works under the assumption of an unlimited defender perception range ultimately fails when the perception is constrained.

2) The addition of ASR: The initial states are as follows: $\mathcal{D}$ is at $\mathbf{x}_D^0 \left( \frac{1}{\sqrt{2}}, -\frac{1}{\sqrt{2}} \right)$ and $\mathcal{A}$ is at $\mathbf{x}_A^0 \left( \frac{3}{2\sqrt{2}}, \frac{3}{2\sqrt{2}} \right)$. The speed ratio is $\nu = 0.47$, and the target radius is 1. In Fig. 5 (a), both sides use the strategy from \cite{r21}, and the attacker successfully breaches the defender's interception to reach the target. However, after introducing an attacker perception constraint (ASR) with a radius of 0.2, as shown in Fig. 5 (b), the defender successfully intercepts the attacker. This suggests that the inclusion of ASR transformed a failed interception into success. The strategy for successful evasion, based on the assumption of the attacker’s unlimited perception range, ultimately fails when the attacker’s perception is limited, preventing timely avoidance of the defender.

Simulation results show that adding perception constraints enhances the model's realism. In practical scenarios, agents often lack complete information and must adapt their decisions, improving their ability to respond in dynamic, uncertain environments. perception constraints also drive agents to develop more complex and flexible strategies, leading to more diverse outcomes and strategic interactions.

\section{Conclusion}
This letter addresses the challenge of patrolling the perimeter of a circular target area for protection, considering perception constraints. We model scenarios where these constraints impact both factions. The game is divided into three stages, each with distinct strategies. When the defender intercepts the attacker during the Full-information Stage, an additional escape stage is introduced for the attacker, providing strategies for both parties.

A promising direction for future research involves exploring scenarios where the attacker can conduct repeated attacks on the target area, using strategic patterns to influence the defender's positioning and gain an advantage for subsequent attacks.

\addtolength{\textheight}{-12cm}   





\end{document}